\documentclass[]{rspublic}

\usepackage[pdftex]{graphicx}

\usepackage{amsthm,amsmath,amssymb,mathrsfs}
\usepackage{bm,bbm}
\usepackage[numbers]{natbib}

\usepackage{palatino}

\def\c2mb{{\rm c}_{2mb}}
\def\s2mb{{\rm s}_{2mb}}
\def\ch2mb{{\rm ch}_{2mb}}
\def\sh2mb{{\rm sh}_{2mb}}

\def\ber{\begin{eqnarray}}
\def\eer{\end{eqnarray}}
\def\be{\begin{equation}}
\def\ee{\end{equation}}
\def\bea{\begin{eqnarray}}
\def\eea{\end{eqnarray}}

\begin{document}

\title[A Geometric Knotspace Template]{A Geometric Knotspace Template}

\author[C.D. Modes, M.O. Magnasco]{Carl D. Modes$^1$ and Marcelo O. Magnasco$^1$}
\affiliation{$^1$Laboratory of Mathematical Physics, Rockefeller University;1230 York Avenue, New York, NY, 10065, USA}



\maketitle

\begin{abstract}{knot theory | medial set | knot complexity | computational knot theory | knotspace deconstruction | knot group}
Early last century witnessed both the complete classification of 2-dimensional manifolds 
and a proof that classification of 4-dimensional manifolds is undecidable,
setting up 3-dimensional manifolds as a central battleground of topology to this day. 
A rather important subset of the 3-manifolds has turned out to be the {\em  knotspaces}, the manifolds left when a thin tube around a knot in 3D space is excised. 
Given a knot diagram it would be desirable to provide as compact a description of its knotspace as feasible; hitherto this has been done by computationally tessellating the knotspace of a given knot into polyhedral complexes using {\em ad hoc} methods of uncontrolled computational complexity.
Here we present an extremely compact representation of the knotspace obtainable directly from a knot diagram; more technically, an explicit, geometrically-inspired polygonal tessellation of a deformation retract of the knotspace of arbitrary knots and links. Our template can be constructed directly from a planar presentation of the knot with $C$ crossings using at most $12C$ polygons bounded by $64 C$ edges, in time $\mathcal{O}(C)$. We show the utility of our template by deriving a novel presentation of the fundamental group, from which we motivate a measure of complexity of the knot diagram. 

\end{abstract}


Few branches of mathematics exemplify as well as the theory of knots the transition from charming pastime to cornerstone importance in entire fields of both the pure and applied domains. The realisation that the knotspaces are a central subset for the classification of 3-manifolds gave knot theory a permanent and central place in topology \cite{Thurston}. The development of knot polynomials tied knot theory to quantum field theories in mathematical physics \cite{Witten} and statistical mechanics \cite{field theory}, and knot theory rapidly became of interest in the theory of polymers \cite{statmech1} and even in wavefront optics \cite{Dennis, protein}. This applied impetus stimulated, in addition to true topological invariants, a lot of activity in variational characteristics, for example the rope length problem \cite{rope length}. 
The discovery that the DNA of living cells becomes entangled, and that life has evolved enzymes, called \textit{topoisomerases}, charged with managing the topology of their DNA \cite{topoisomerase} gave an applied and urgent impetus in a completely unforeseen area: one of the most important classes of antibiotics today, the \textit{ gyrase inhibitors}, act by tampering with bacterial enzymes (gyrases) which control the topology of bacterial DNA.

Given this interest, the computational complexity of knot theory has lagged development. Haken conceived an algorithm that would classify all knotspaces by cutting them along special surfaces, called {\em incompressible}, preserving special markings that would permit their re-gluing \cite{Haken}, a tour-de-force program advanced by Hemion \cite{classification1} and further \cite{classification2, classification3, tabulation, spines}. The extant proof shows the algorithm terminates in finite time; it has no proven running time, and in fact as far as we are aware has never been fully implemented. As an example, the \textit{unknot recognition} problem seeks to deduce, from a planar presentation of a knot, whether it is genuinely knotted or merely a convoluted presentation of the \textit{unknot}. The precise complexity status of this problem is still unsettled; use of the invariant knot polynomials does not avail, as they do not guarantee that no knot has nontrivial polynomial and run in exponential time anyway. A polynomial-time algorithm has been claimed using a reduction to braids \cite{braids}.
The theory being thus unsettled, it is unsurprising that the current state of the art in computational analysis of a knot- or linkspace is complicated and \textit{ad hoc} \cite{field theory, stateart2}, 
leading to issues in generality and scaleability of these analyses.  

We introduce as a solution to these issues a geometrically derived, algorithmically simple, systematic deconstruction, or ``template,'' of a knot- or linkspace.  Our template is homotopic to the {\em medial set} of the knot, a set of 2D surfaces in 3D which, in a sense, separate the knot from itself.  We shall prove this medial set to be a {\em deformation retract} of the knotspace.  Our template is an explicit computation of this medial set from the customary planar presentation of the knot or ``knot diagram'', by assuming the ``depth'' of the figure to be much smaller than its lateral dimensions, and permitting minimal inessential deformations of the surface to smooth it onto extrinsically-curved polygons. While this construction is virtually a pencil-and-paper one, it is well-suited for algorithmic implementation and runs extremely efficiently. 

We use our template to derive a novel presentation of the fundamental group of the knot, and then use it to define a measure of complexity of the presentation.
The fundamental group of the knotspace is the group of closed loops in the knotspace starting and ending at a chosen base point (the choice is immaterial). A classical construction, the \textit{Wirtinger presentation} of the fundamental group \cite{Wirtinger}, is derived directly from the knot diagram by dividing the drawing of the knot diagram into strands, which are considered to terminate whenever they pass under another strand in the knot diagram.  This group presentation has as many generators as strands, i.e. $C$ generators, and $C$ relationships having exactly $4$ letters each. At least 1 relationship has to be dependent on the others, as the simplest possible fundamental group that can be associated to a knotspace is the free group on a single generator (in the case of an unknot). 

As we shall show, an application of the van Kampen theorem allows us to use our template to derive a presentation of the fundamental group that is, to some extent, a two-sided version of Wirtinger. Our presentation has $2C$ generators (or, `letters') with $2C$ relationships. However, the number of letters in each relationship \textit{varies}, with an average taken over the set of $2C$ of between $3$ and $4$ letters. When a relationship has a single letter, it has the form $g=1$, and the corresponding generator immediately cancels.  When a relationship has two letters, it either has the form $g=h$ or $g=h^{-1}$.  In either case, one generator may be used to replace the other and again a generator vanishes from the algebra.  Computationally the worst possible scenario is thus when all relationships have exactly 3 or 4 letters, satisfying the average without giving rise to simple cancellations or replacements: this happens only when the knot is alternating. We shall show below that we can define a measurement of the {\em presentation complexity} by taking the \textit{geometric mean} of the relationship lengths--when this geometric mean is low, then the spread in presentation lengths gives rise to a number of short relations that cause the presentation to simplify. The presentation complexity is therefore maximal for alternating presentations. 

Furthermore our use of van Kampen will show that our presentation of the fundamental group is flexibly determined by the details of a pair of quotient graphs, and the way that the generators for those quotient graphs' simpler fundamental groups sit inside the algebra of the medial graph's fundamental group.  This flexibility allows for the generation of a spectrum of equivalent algebraic relations and may thus aid in simplification and identification algorithms. 

We provide for the general reader a brief practical summary of our procedure and the multiply branched surface that results from it in Box 1, and the group presentation in Box 2, complete with worked examples.  In the remaining body of the paper, we delve deeper into the motivations, derivations, examples, and implications that accompany the template.

\begin{table}
{\bf BOX 1: Summary of templating algorithm with the trefoil as example}\\
\begin{tabular}{| l p{8cm}  |}
\hline
\raisebox{-3cm}{\includegraphics[width=4cm]{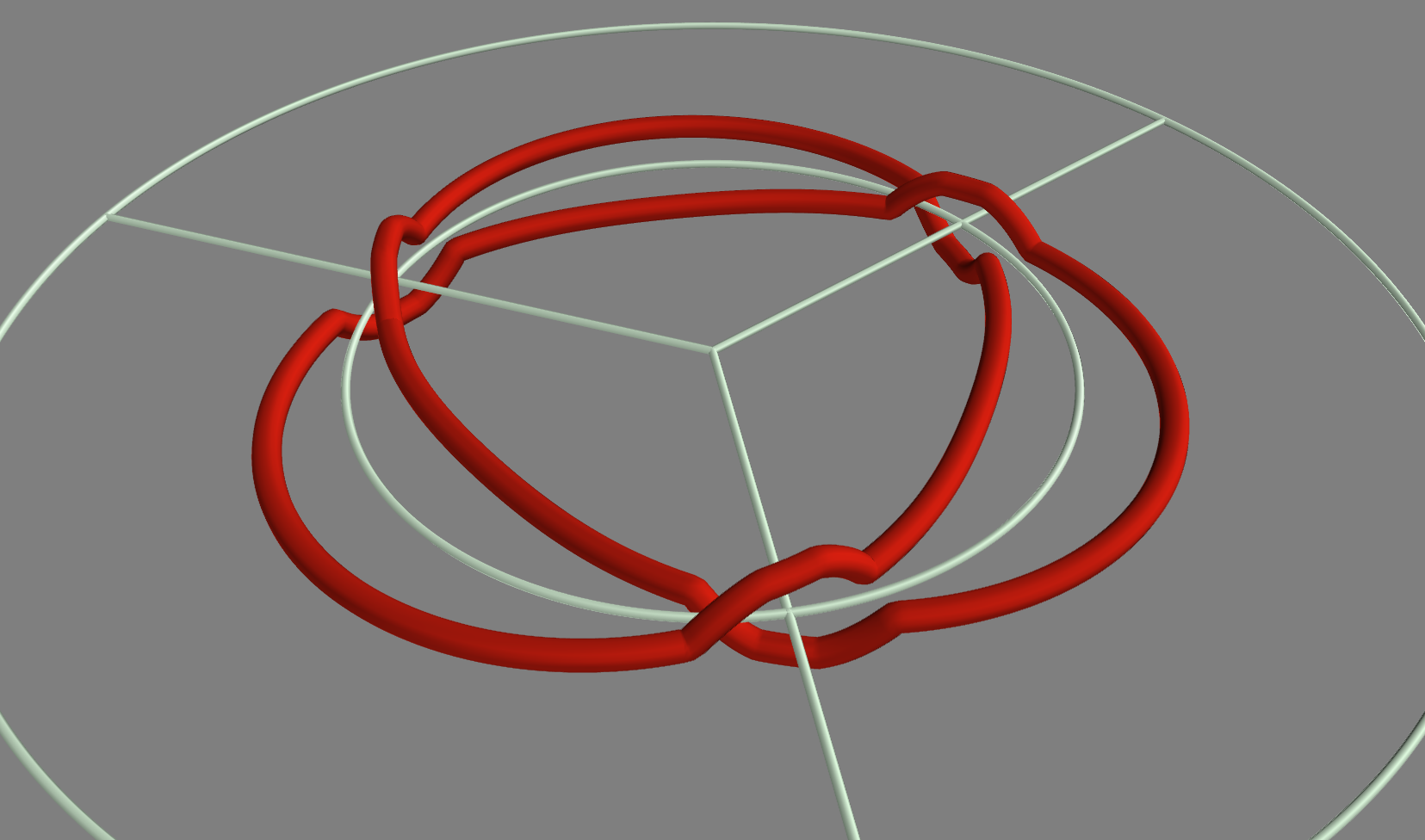}}&
 
 We start from an almost flat presentation of the knot in the x-y plane (red); it will have some given number of crossings $C$, of which $E$ lie at its exterior boundary and $C-E$ lie inside (Here $C=3$, $E=3$). We then construct the \textit{ bounded medial graph} (pale green) of the knot., made by placing an outer bounding circle on the medial graph, on which the unbounded segments of the medial graph terminate.  The medial graph itself is simply the pruned deformation retract of an almost flat knot diagram where the crossings have been replaced by four-fold vertices in the plane. \\
\hline

\raisebox{-3.7cm}{\includegraphics[width=4cm]{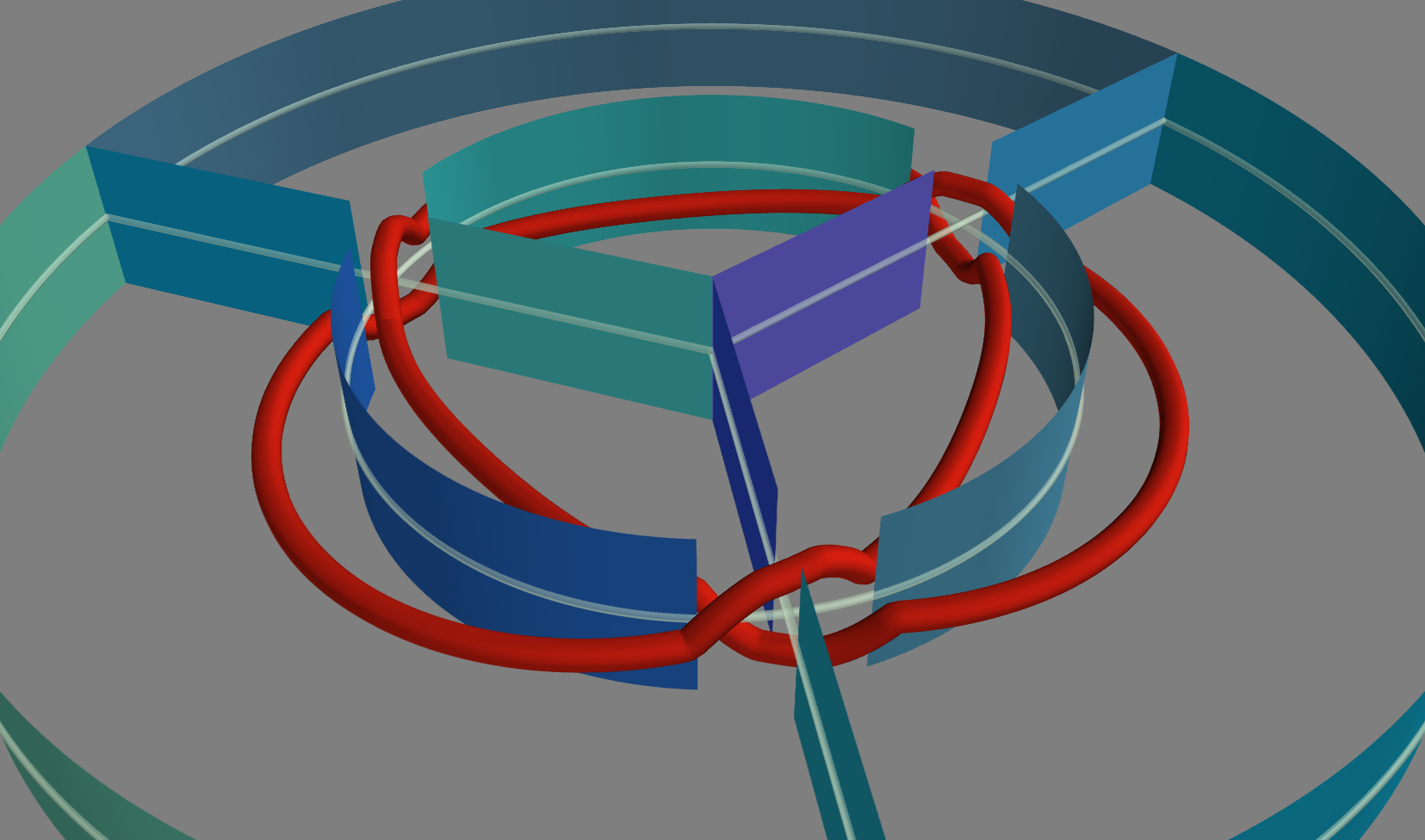}}&

{\bf Step 1: Medial walls}. We exclude from the bounded medial graph small neighbourhoods around the knot crossings, and ``extrude'' vertical walls out of the remainder. Every wall at this stage is topologically a rectangle, which may be curved in the $(x,y)$ plane but straight vertically.  Each facet $F$ of the knot diagram contains a star-like arrangement of as many walls as crossings around the facet (the central facet here has three walls impinging), unless the number of crossings is 2, in which case the two segments can be collapsed into a single segment (the three outer facets have a single curved wall each).  
There are $4 C -T $ internal medial walls, where $T$ is the number of two-sided facets of the knot diagram, plus $E$ walls into which the outer ring is subdivided wherever the outer medial walls impinge upon it. \\
\hline

\raisebox{-2.4cm}{\includegraphics[width=4cm]{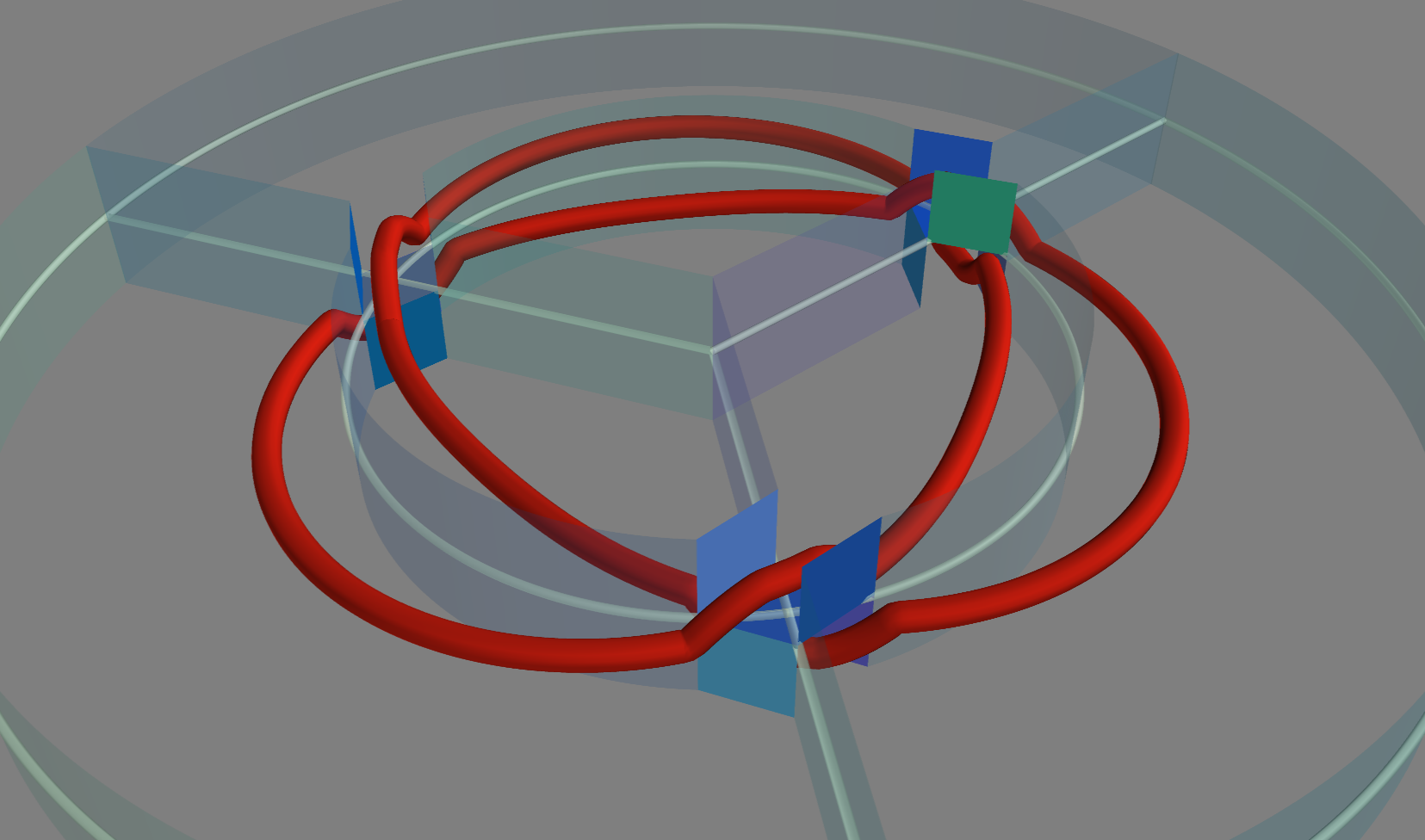}}&

{\bf Step 2: Saddles}. A 5-rectangle ``saddle'', formed from a central square with four rectangles attached to each edge, is fitted into each knot crossing.  The square lies in the plane of the knot, ``in" the crossing, with each vertex over a point of the medial graph. The long rectangles are folded perpendicular to this plane, alternating up and down, so that the saddle ``cradles'' the crossing knot segments. There are evidently $5 C$ rectangles involved. \\
\hline
\raisebox{-4cm}{\includegraphics[width=4cm]{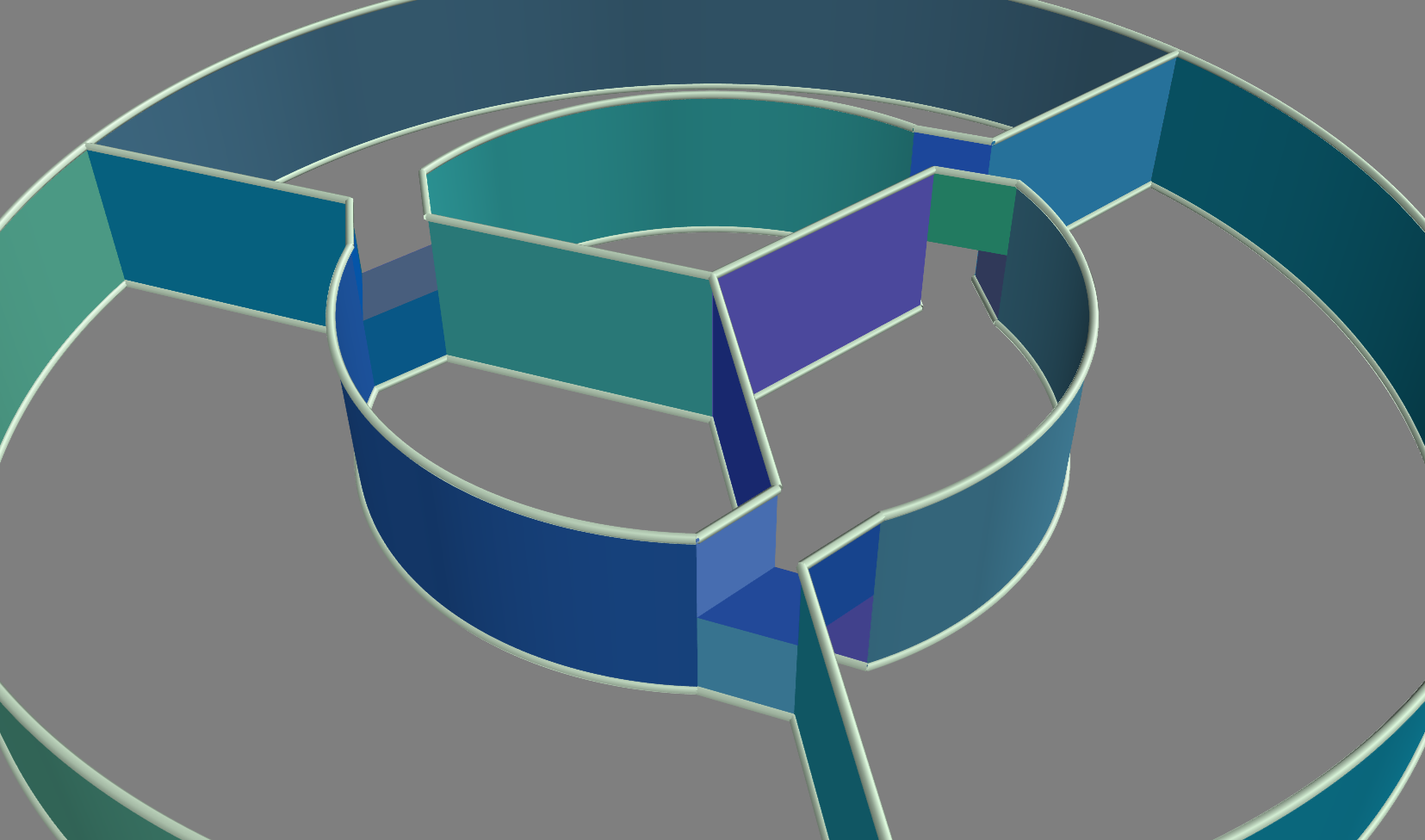}}&

{\bf Step 3: Lids.} Top and bottom circular ``lids'' are placed to close the construct, i.e. two disks. The top and bottom lids are subdivided into facets where the medial walls and the saddles impinge on them; we call these intersections the \textit{ upper and lower graphs}, outlined above in light green. This step adds a total of $2 C$ facets, each having an average of $11$ edges. This completes our construction. Collecting facet counts from each section of the template yields at most $11C + E$, allowing for simple ``bend" curvature of some of the faces, of which $9C+E$ are four-sided, while $2C$ are on average $11$-sided.  This collection of polygons together with the information of how they impinge upon one another contains all the topological information of the knotspace, as a minimal, infinitesimal thickening of each polygon into thin polyhedra recovers a complete polyhedral decomposition of the knotspace.
\\
\hline
\end{tabular}
\end{table}

\section*{Construction of the template}\label{sect:template}

The fundamental conceit of this paper is 
to consider a geometrically-informed deformation retract of the knot space that allows for a complete collapse onto a set of joined surfaces which we call the \textit{knotspace template} or simply, the template.  Since deformation retracts preserve the homotopy of paths \cite{Hatcher}, the fundamental group in particular will be a topological invariant of this transformation.  On the other hand, homological spectra, relying as they do on cell-complex constructions, will not be preserved under the dimensional collapse, and as a result the template will not be \textit{strictly} sensitive to some popular knot invariants that require them, such as Floer Homology \cite{Floer1, Floer2}.  However, a simple ``re-inflation" of the polygons of the template into thin rectangular prisms restores the missing dimension in a topology-preserving way and homology may proceed from there. 

We shall proceed through the following stages. First, we will show that the medial set of the knot, which can be evaluated through Voronoi procedures, is a proper deformation retract. Second, we shall specialize the construction to a quasi-two-dimensional case, in which we examine a knot diagram having an infinitesimal depth: we shall show that away from the crossings, the medial set coincides with the medial set of the perfectly flat knot graph, and then show that in a small neighborhood of the crossings the medial set is a saddle-like surface cradling the strands, which we can deform into a polygonal saddle. Putting everything together requires a few extra surfaces to compactify the construction, in the shape of an outer ring and a top and a bottom lid. 

\subsection*{Voronoi tilings and the medial set}\ \\
In order to specify the specific deformation retract onto the template we must be able to prescribe a flow from every point of the knot space to a (possibly branched and faceted) surface immersed in it.  This surface will be the deformation retract.

Voronoi tilings of spaces, a well-worn tool familiar in statistical mechanics \cite{Kirkwood, stat} and the mathematics of packings \cite{Voronoi}, offers a conceptual way forward.  Imagine that the knotted curve whose knotspace we wish to deformation retract is replaced by an evenly spaced set of points.  Note that the Voronoi tiling associated with these points will consist of two distinct classes of region boundaries -- those arising from adjacent points along the path of the knot and those arising from non-adjacent points.  Clearly the boundaries arising from adjacent points do not include information inherent to the knotted curve -- as they are generated locally -- and may be discarded.  The remaining boundaries comprise a nascent form of the template.  Unfortunately, choosing the correct boundaries is algorithmically difficult and indeed, the status of the mapping as a deformation retract is cast into question if there are several extraneous target surfaces that must be ignored.  What is to be done?

Happily, there is a generalization of the Voronoi procedure to higher dimensional objects that obviates the need to reduce the knot to a set of points in the first place.  The medial set of geometric analysis \cite{medial} or the spine of manifold recognition algorithms \cite{spines} is this generalization.  Such a set is defined in an ambient space of dimension $n+1$ as the locus of $n$-ball centers for the set of $n$-balls that are tangent to the hypersurface of interest -- the knot in our case -- in at least two distinct places and do not otherwise internally intersect it.  For complicated curves, the medial set may develop many small spikes and protrusions, but these do not change the topological character of the retract and may be safely pruned.

Certainly in the case of distinct points, the procedure for generating the medial set recapitulates the Voronoi boundaries, as multiple tangencies ensures the $n$-ball center is between Voronoi regions and the condition against internal intersection ensures that there is no other closer point's Voronoi region to which the $n$-ball center should belong.  The specification of the deformation retract map from a closed curve to a surface is straightforward: any point that lies anywhere along a radius associated with a point of multiple tangency is mapped to the center of the corresponding sphere.  All other points are mapped to the surface at infinity.  This map is well-defined as a point cannot lie on multiple such radii and continuous as the family of multiply tangent balls can be continuously traversed by inflation or ``rolling" on the curve.  The requirement that points already in the image are stationary is satisfied by construction.

\begin{figure}[!t]
\ \\
\includegraphics[width=4.1cm]{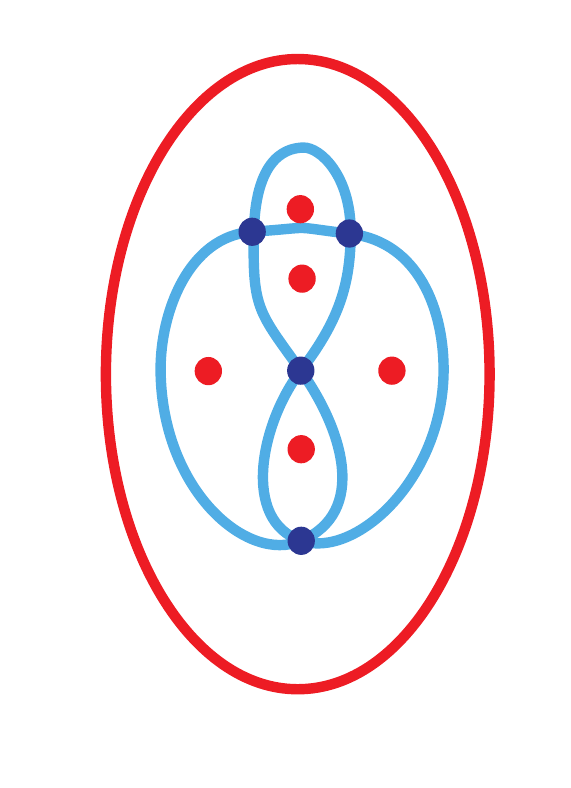}
\includegraphics[width=4.1cm]{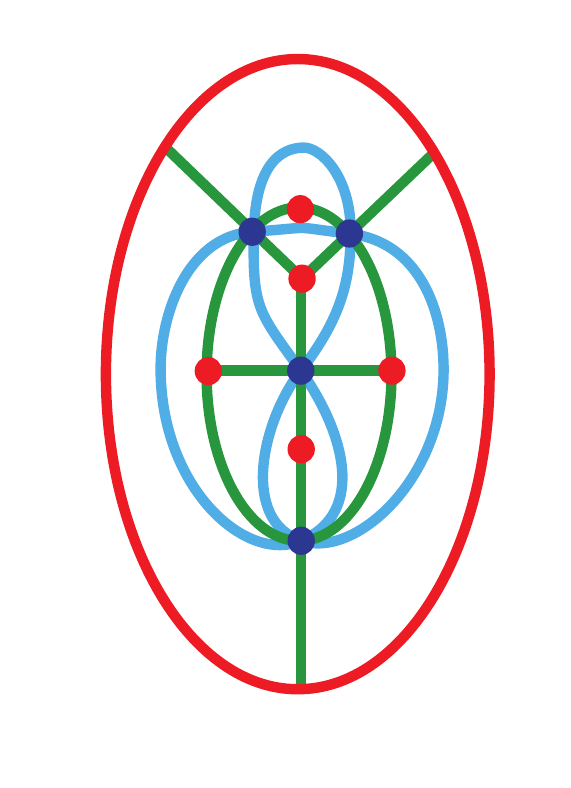}
\includegraphics[width=4.1cm]{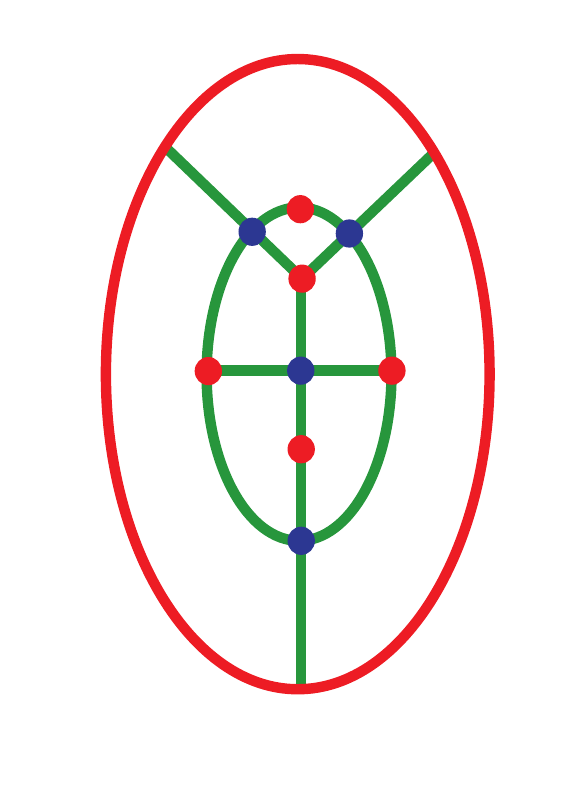}
\caption{The construction of the bounded medial graph, illustrated for the {\em figure eight} knot. Left, the knot's diagram is projected onto a planar graph of connectivity four, its fourfold vertices highlighted in dark blue. Dual vertices (our ``stars'', shown in red) are added for every facet of the original graph. The outer unbounded facet gets a vertex which is here shown unraveled as a surrounding red circle. Center, the blue fourfolds are connected to the red stars whenever the facet of the original graph has that fourfold as vertex, shown in green. This is the full bounded medial graph: a bipartite graph, one of whose partitions are the original fourfold vertices, the other partition corresponding to the facets. Each facet of the medial graph contains a piece of a knot strand. Right, the medial graph by itself. }
\label{fig:medial}
\end{figure} 

\subsection{The bounded medial graph}\ \\
A {\em knot} or {\em link diagram} consists of a planar projection of the knot or link, in which the choice of overpass/underpass is shown diagrammatically as a continuous overpass and an interrupted underpass strand, creating the illusion of an occlusion in 3D. Abstractly, the knot diagram is therefore a planar graph of connectivity four, equipped with a (Boolean) choice of overpassing strand at each vertex. (Technically, it is often a multigraph, meaning that two nodes may be connected by multiple strands). The number of components in the link or knot is given by the 4-fold graph regardless of the choices of overpass. Evidently there are infinitely many knot diagrams for each knot. 

The medial graph of the knot diagram is constructed as a graph abstraction of the pruned medial set; after inessential pruning and simplification it can be cast into the specific form we shall henceforth use, illustrated for the figure eight knot in Figure 1. The medial graph is a {\em bipartite} planar graph; the vertices of one partition are 4-fold vertices identified 1-to-1 with the original vertices of the knot diagram, referred henceforth as the {\em fourfolds}. The four edges emanating from these vertices are not identified with the knot strands but rather with the quadrants between them, and hence appear at $45^\circ$ from the original edges. The second set of vertices corresponds 1-1 with the facets of the original knot diagram; if the facet is a $p$-sided polygon then the corresponding vertex is a $p$-fold star, and each of its $p$ edges connects with each of the $p$ four-fold vertices surrounding the facet. We call the members of this second partition the {\em stars}. 

Therefore the fourfolds of the medial graph correspond to the original fourfolds of the knot diagram, the stars correspond to the facets of the knot diagram, and the facets of the medial graph correspond to the links (or strand segments) of the knot diagram. It necessarily follows that each link of the medial graph, being the boundary between two facets, corresponds to a way to separate two knot strands from each other. 

The fourfolds that lie on the outer boundary of the diagram have one edge emanating towards the outside, in principle to infinity. We compactify the diagram by surrounding it with a circle representing a star at infinity, on which these outer unbounded edges now terminate. We will depict the star representing the unbounded outer facet of the knot diagram not as a point but by unraveling it as a circle, though we keep in mind it is a vertex to fully preserve the bipartite structure of the medial graph. 

In principle it is possible to prune this construction further in the case of facets bounded by 1 or 2 edges; in the first case the vertex of connectivity 1 and its edge can be removed; in the second case, the edge of connectivity 2 can be removed and the edges impinging on it merged. However, this further pruning destroys the bipartite structure of the medial graph. 

Not every bipartite planar graph, one of whose partitions consists exclusively of 4-fold vertices, is a valid medial graph of a knot or link diagram. A necessary and sufficient condition is that every facet of the bipartite graph be 4-sided. It is necessary because the construction of the medial graph makes them so: since every strand segment corresponds to a facet, and strand segments terminate in two fourfolds, facets naturally have two fourfolds and two stars as their boundary vertices. It is sufficient because when a facet is 4-sided, by the alternating property it necessarily has two fourfolds. Connecting those fourfolds gives back the knot diagram. (We have glossed over the inconvenient if irrelevant case of a fourfold connected to itself in the knot diagram, removable by a simple Reidemeister move). Therefore, the set of all valid medial graphs can be expressed in closed form as a bipartite planar graph, one of whose partitions consists of fourfold vertices, and whose facets are all four sided. 

Using both the medial graph {\em and} the overpass choices of the knot diagram, we can define two graphs called the {\em upper} and {\em lower} graphs. Their vertices are the stars of the medial graph; their edges are created by removing the fourfolds, connecting together the edges impinging on the fourfold in the usual manner employed in skein diagrams; more specifically, the two pairs of adjacent edges associated with the underpass are connected for the {\em upper} graph, and the pairs associated with the overpass are connected for the {\em lower} graph. These two graphs will appear naturally as an important part of our template; as shown in Box 1, Step 3 these graphs trace the intersection of the template with the top and bottom lids. 


\subsection*{Knot crossings and saddles}\label{subsect:saddle}\ \\

The medial set of a quasi-flat knot coincides, away from the crossings, with the medial set of the flat knot graph; it is therefore a set of surfaces erected perpendicularly to the medial graph. The medial graph contains by construction no information about overcross/undercross choice. But in a small neighborhood of every crossing we need to consider the full impact of the strands in 3D and the choice of overpass.  From the previous discussion we may expect that whatever happens with the medial set at these crossings must affect a ``hooking together" of the medial-set planes that result from the knot regions far from a crossing. We shall now show that the medial set in a neighborhood of a crossing is simply a saddle, which indeed hooks together the pieces of surface already erected and contains the topological information relevant to encoding the overpass/underpass choice at every crossing. 

\begin{figure}[!t]
\ \\
\includegraphics[width=6.2cm]{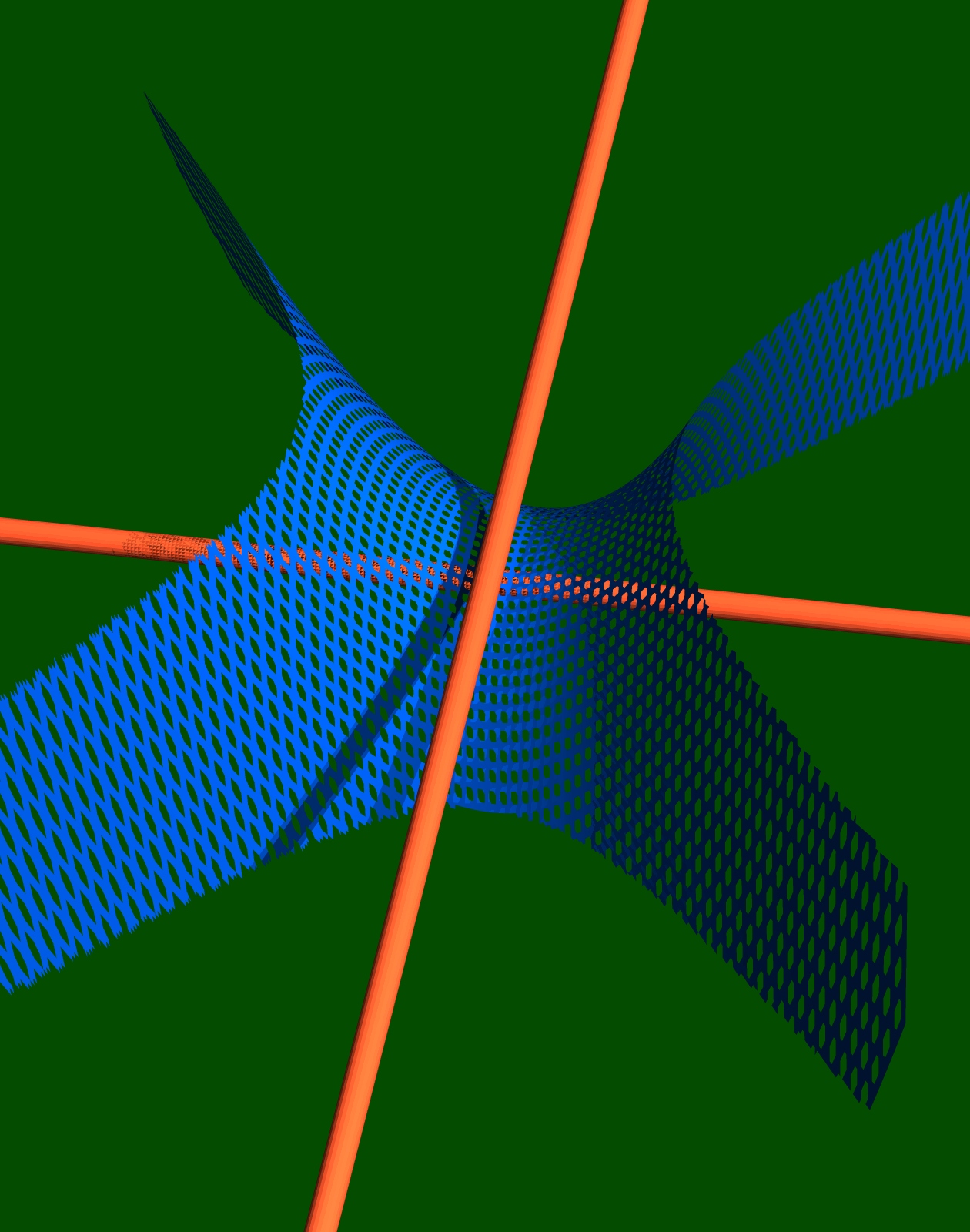}
\includegraphics[width=6.2cm]{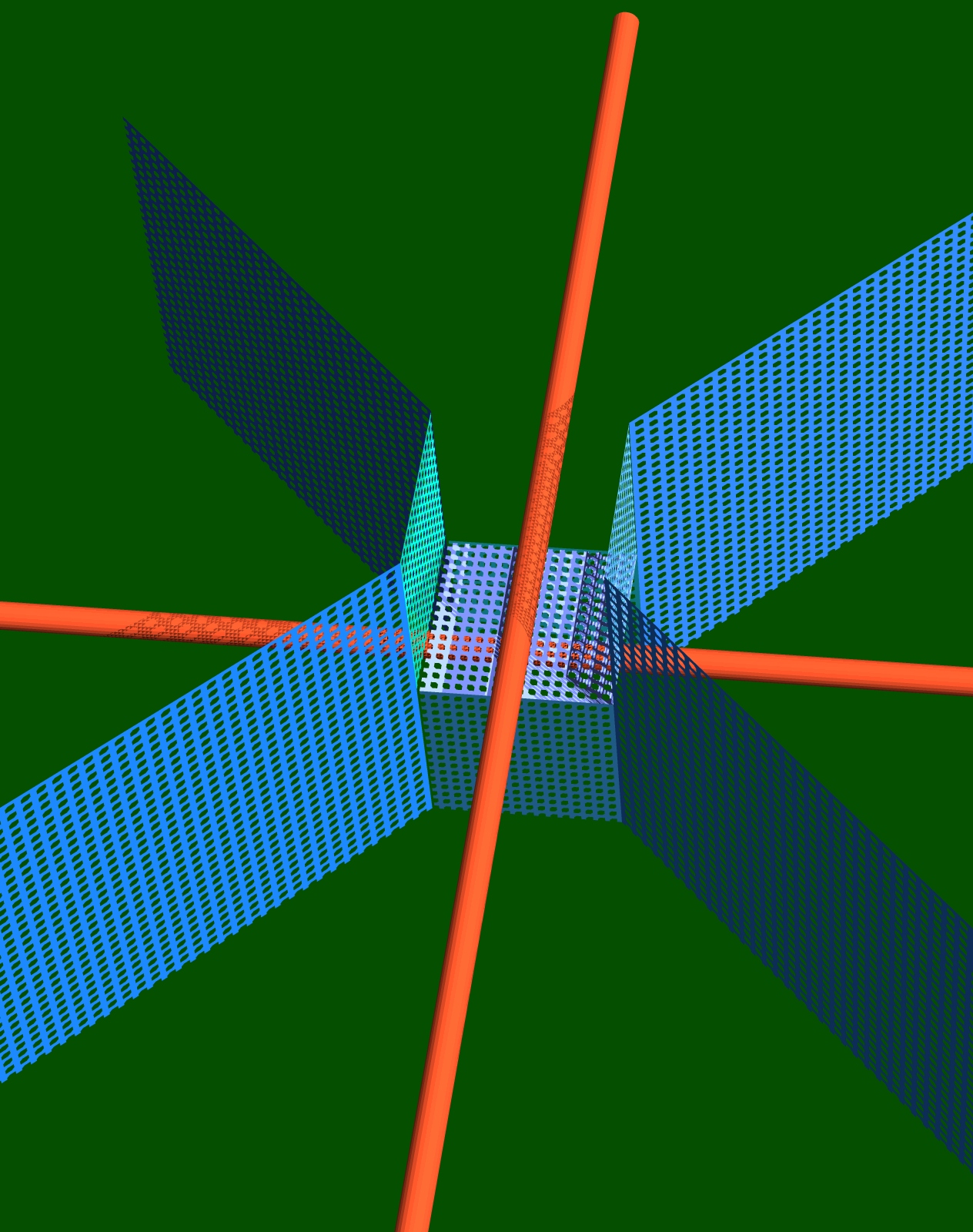}
\caption{The generation of the medial set near a crossing in a quasi-planar knot diagram (left) and subsequent tesselation into simple rectangular pieces (right, see also Box 1).  The initial smooth surface is locally a hyperboloid.}
\label{fig:sheath}
\end{figure} 

Let us model such a crossing as a pair of skew lines separated at their closest passing by a small distance, $\delta$.  In such a scenario, one may find the medial set by imagining a pair of cylindrical sheathes with arbitrarily small radius centered, respectively, on each of the skew lines.  These sheathes' radii are then continuously inflated at the same rate and the locus of points over which the two sheathes intersect is the medial set.  The resulting surface is a hyperboloid that cradles the pair of generating skew lines (Fig. 2).  Note  that far from the crossing along a direction in the thin plane containing the lines, this surface twists around to form a set of four vertical half-planes at an angle of $\pi/4$ to the crossed lines. These half-planes may then seamlessly connect to the polygonally-generated branches of the medial set discussed above, as required.  Meanwhile, very far from the crossing in the direction normal to the thin plane containing the lines the surface manifests as two separate vertical sheets folded through an angle of $\pi/2$ over a distance comparable to $\delta$.  These folded sheets are closest to one another at their respective folds, and are separated by a distance, again, comparable to $\delta$.  Which pairing of the four half-planes that are connected by a fold switches above and below the plane of the crossing, correlating with which of the two perpendicular cylinders is farther from the midpoint of the crossing.

\subsection*{Examples of the template}\ \\
Given a quasi-planar presentation of a knot or unknot the template will always be composed of four distinct regions, examples of which are presented in Box 1.  The first of these, above the plane of the quasi-planar presentation, consists of vertical sheets in the manner of the planar curves with four-fold junctions discussed previously.  In this case, however, the sheets do not emanate symmetrically from the junction (now knot crossing) -- rather, each pair of sheets on the same side of the upper section of the knot crossing are smoothly joined together and the pairs of sheets across opposite sides of the upper section of the knot crossing are not connected at all.  The second region of the template is similar to the first, but below the plane of the presentation.  Now the vertical sheets are connected in the opposite manner, with the pairs on the same side of the lower section of the knot crossing smoothly joined together and the pairs separated by the lower curve not connected at all.  The third region of the template is the region containing the neighborhood of the plane of the quasi-planar presentation -- where the vertical sheets of the first two regions must be joined together.  Again, away from the knot crossings in this region the template is simply the vertical sheets of the medial sets, as it was in the previous regions.  However, near the knot crossings, the template is saddle shaped, as depicted in Fig. 2.  By adopting this shape, the template both correctly captures the region between the crossed strands of the knot and connects to the sheets both above and below the plane of the knot diagram in a way that is consistent with the way the vertical sheets connect or do not connect in the first two regions of the template.  Finally, the fourth region of the template is the surface at infinity, which may be imagined to be a sphere, or, more usefully, a cylinder or hat box.  Where the sheets of the first two regions of the template meet the surface at infinity there is a triple junction.

In addition to the template of the trefoil as shown in Box 1, many other simple examples are easy to produce, such as a trefoil-like unknot, the figure eight knot, a figure-eight-like unknot, and even links, such as the Hopf link or the Borromean Rings (Fig. 3)

\begin{figure}[t]
\begin{center}
\resizebox{9cm}{!}{\includegraphics{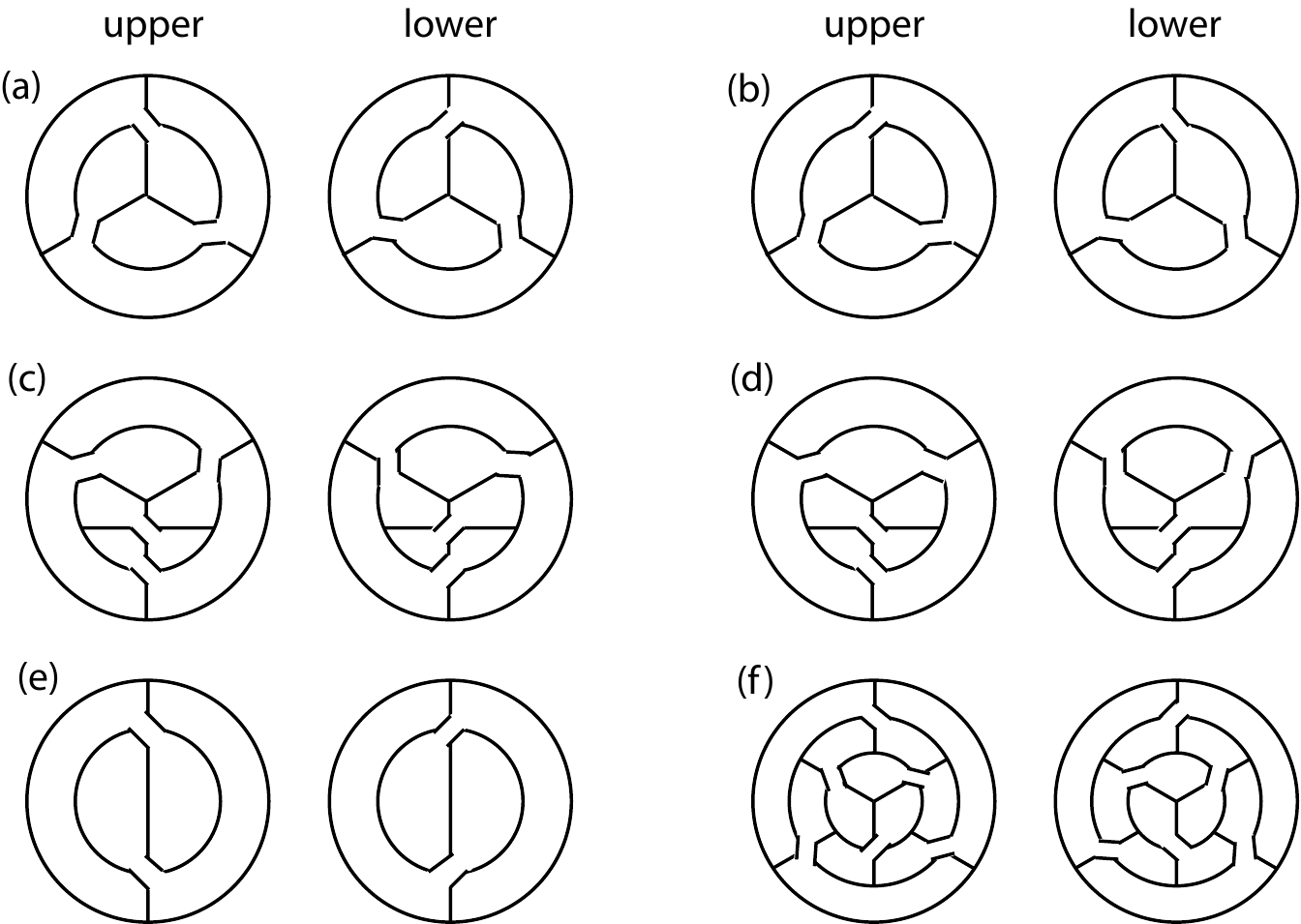}}
\end{center}
\caption{Several examples of the template for different knots, presentations of the unknot, or links.  Each example is depicted as two slices parallel to the plane of the knot diagram, one through the upper region of the template, the other through the lower region.  (a) The trefoil knot.  (b) A trefoil-like presentation of the unknot.  (c) The figure eight knot.  (d) A figure-eight-like presentation of the unknot.  (e) The Hopf link.  (f) The Borromean Rings.}
\label{fig:KTexamples}
\end{figure}

\section*{The fundamental group}\label{sect:fund}

We now demonstrate that the template provides a useful route to calculating topological information.  The fundamental group of the knot space, $\pi_1 (K)$, is an obvious candidate for investigation, as the deformation retract procedure that collapsed the knot space explicitly preserves homotopies.  Though $\pi_1 (K)$ is a useful topological invariant, it alone cannot fully differentiate knots, as, for example, the knot sum of a square knot and a granny knot has the same fundamental group as the knot sum of a square knot and the mirror of the granny knot, despite these two knot sums producing distinct knots.  On the other hand, $\pi_1 (K)$ does distinguish prime knots \cite{pi1prime}, and in particular, will distinguish the unknot from a knot.

In order to proceed we must separate the template into more manageable chunks.  The venerable van Kampen theorem \cite{Hatcher} is the appropriate tool for this job, relating the fundamental groups of two open sets whose intersection is connected and whose union is the space of interest.

We settle on three open covering sets (Fig. 4) -- the upper cap set, $C_U$, covering the top cylindrical cap and all the surface at infinity and vertical planes extending down from the cap almost, but not all the way, to the plane of the knot; the lower cap set, $C_L$, covering the bottom cylindrical cap and all the surface at infinity and vertical planes extending up from the cap almost to the plane of the knot; the mid-band, $B$ encompassing the plane of the knot and all the saddle cross-overs and a section of the vertical planes and surface at infinity extending both up and down, but not reaching the cylindrical caps.  For such a choice of covering sets, the intersections are well characterized -- in fact, we have already done so -- $C_U \cap B$ is a set of vertical planes, open on both vertical ends, whose horizontal cross section is the ``upper" slice of Fig. 3.  Likewise, $B \cap C_L$ is a set of open-ended vertical planes with horizontal cross section given by the ``lower" slice of Fig. 3.

We must also make a choice about the pairwise ordering of the sets to be addressed: we will first calculate $\pi_1(C_U \cup B)$ and then $\pi_1((C_U \cup B) \cup C_L)$.  Hence:

\bea
\pi_1(K) \cong \pi_1(\textit{template}) &\cong& \pi_1(C_U \cup B) \ast \pi_1(C_L) / [\cdot] \label{vK1} \\
\pi_1(C_U \cup B) &\cong& \pi_1(C_U) \ast \pi_1(B) / [\cdot]. \label{vK2}
\eea
where $[\cdot]$ here is $[i_{FG}(\omega) i_{GF}(\omega)^{-1}]$ for $i_{FG} : \pi_1(F \cap G) \rightarrow \pi_1(F)$ the map induced by the inclusion $F \cap G \subset F$ and $F,G$ the two covering sets under that particular application of van Kampen.  The square brackets denote conjugacy under homotopy.

We may leverage the same deformation retract procedure that produced the template to make drastic simplifications of the covering sets under consideration.  First, note that the two capping sets, $C_U$ and $C_L$ may be retracted all the way down to a single point, as the vertical sheets may be smoothly retracted back into the cap, and once the cap is all that remains it may be retracted to a point as well.  Therefore, we have triviality of the associated fundamental groups:

\be
\pi_1(C_U) \cong \pi_1(C_L) \cong 1.
\ee
Furthermore, the intersection bands are more than just characterized by the slices of Fig. 3 -- they deformation retract onto them as well.  Since those slices are graphs, we immediately have:

\bea
\pi_1(C_U \cap B) &\cong& \mathbf{Z} \ast \cdots \ast \mathbf{Z} \\
\pi_1(B \cap C_L) &\cong&  \mathbf{Z} \ast \cdots \ast \mathbf{Z}
\eea
where the number of copies of $\mathbf{Z}$ in each fundamental group is given by the number of independent loops in the associated graph \cite{Hatcher}.  We refer to these graphs as the \textit{upper} and \textit{lower graphs}, respectively.  Since these graphs differ only by the connections over/under the knot crossings, the number of these independent loops will be the same for each, and thus $\pi_1(C_U \cap B) \cong \pi_1(B \cap C_L)$.

\begin{table}
{\bf Box 2. Presentation of the fundamental group}\\
\begin{tabular}{ | r | r |  p{ 0.5 \columnwidth } |}
\hline 
Trefoil & Unknot & \ \\\hline
\raisebox{-0.15 \columnwidth}{\includegraphics[width=0.2 \columnwidth]{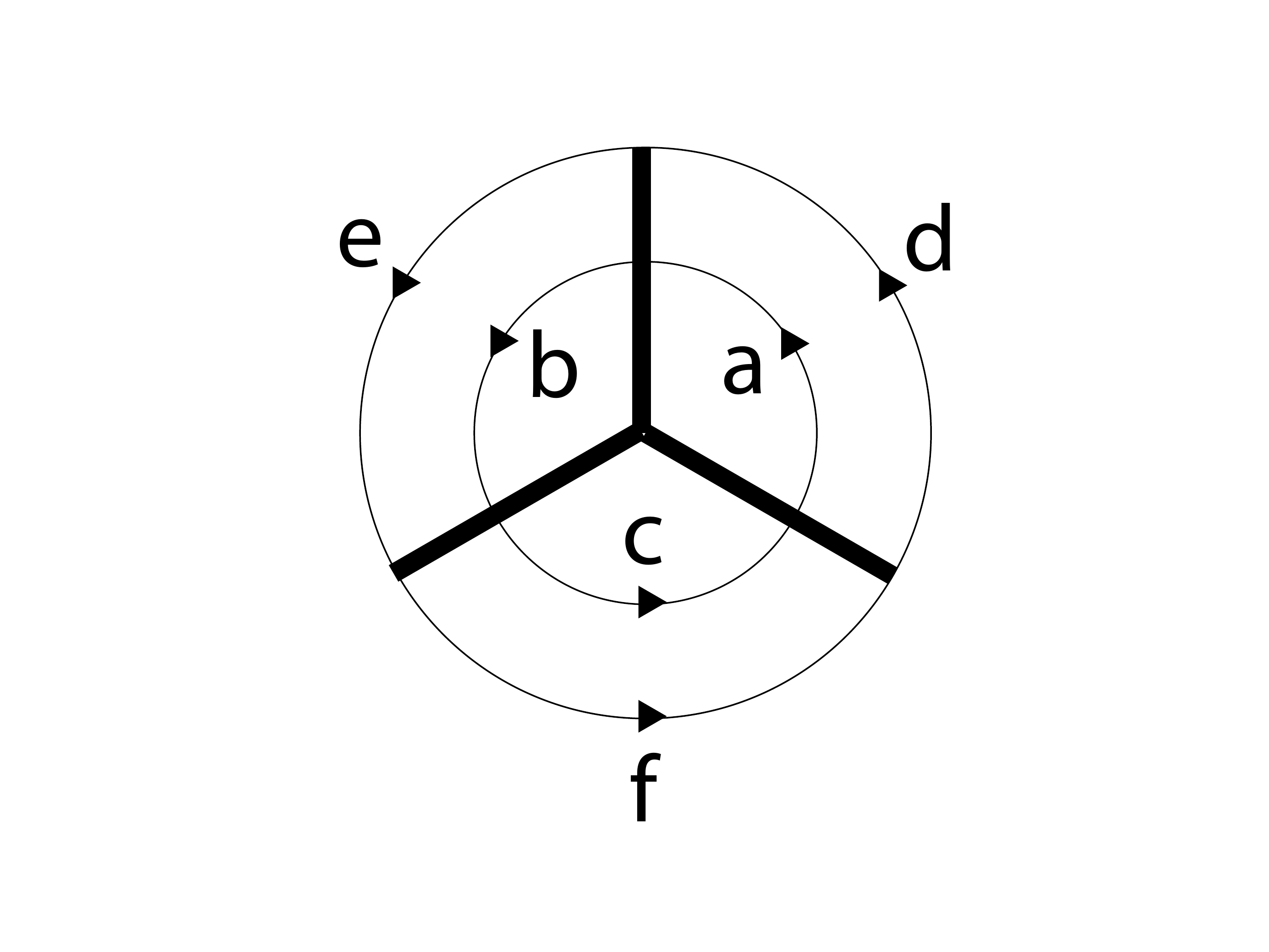}} &
\raisebox{-0.15 \columnwidth}{\includegraphics[width=0.2 \columnwidth]{tabmedial.pdf}} & Step 1: Generate a spanning tree of the medial graph, indicated as bold lines in the figures to the left. All edges of the medial graph not in the spanning tree will be the generators of the fundamental group of the medial graph: orient them and label them. (In this figure we shall orient everything counterclockwise)\\\hline
\raisebox{-0.15 \columnwidth}{\includegraphics[width=0.2 \columnwidth]{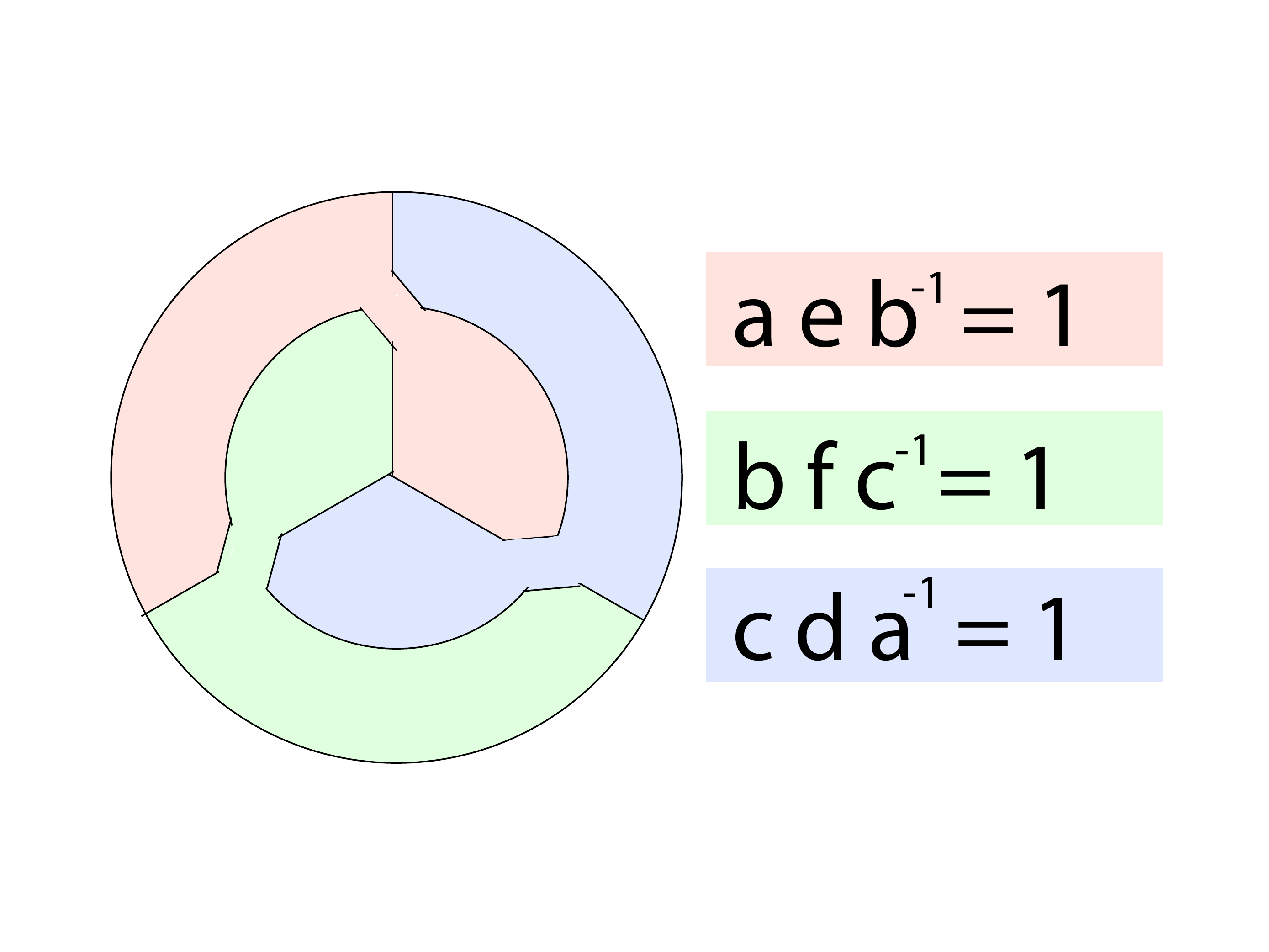}} &
\raisebox{-0.15 \columnwidth}{\includegraphics[width=0.2 \columnwidth]{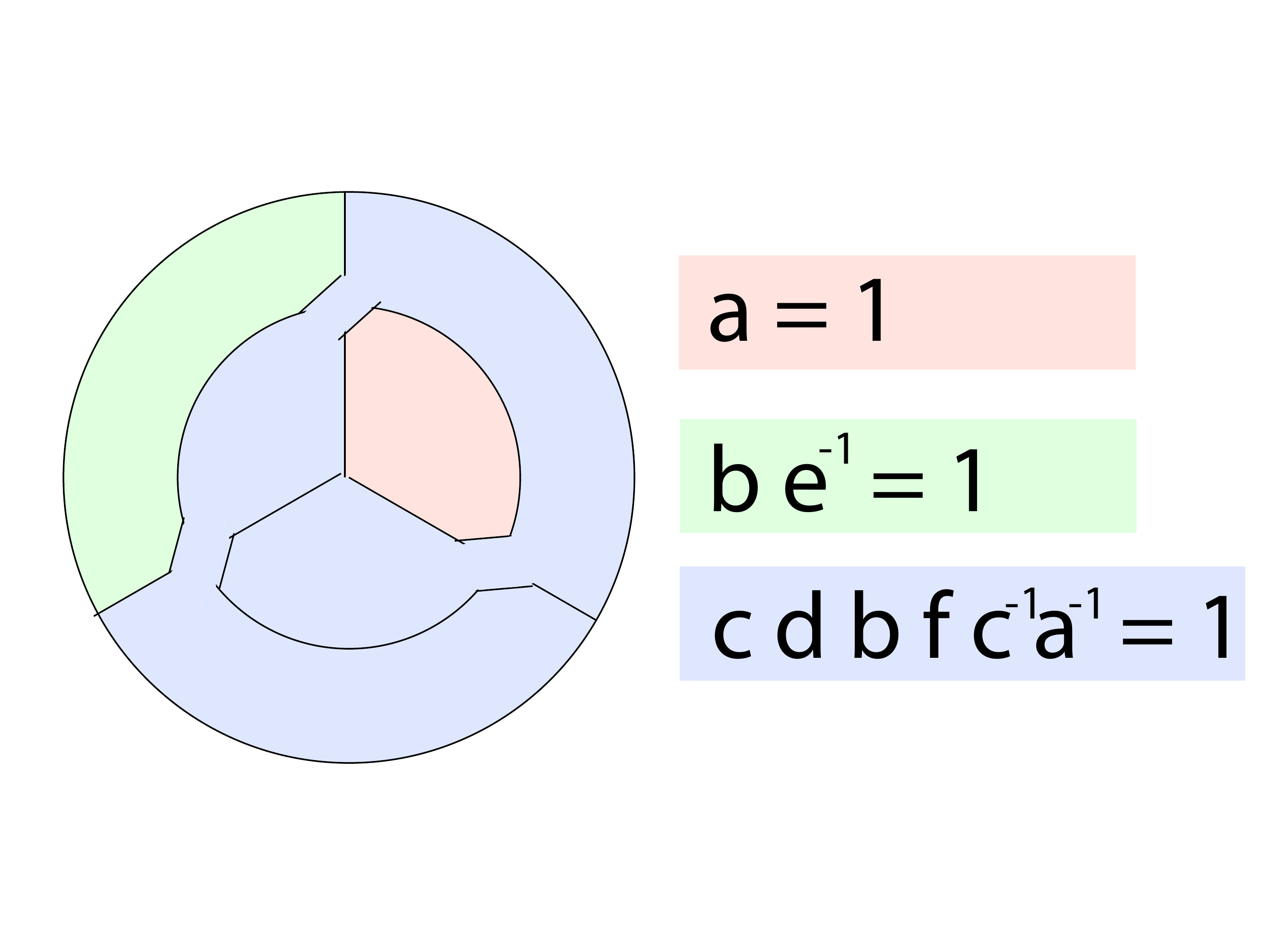}} & Step 2: Every facet in the upper graph gives rise to a relation. Traverse each facet, and record the generators you encounter, raised to the $-1$ power if traversed in the opposite orientation; equate the result to $1$ to express that the free group of the medial graph has been divided by this loop in the upper graph. \\\hline
\raisebox{-0.18 \columnwidth}{\includegraphics[width=0.2 \columnwidth]{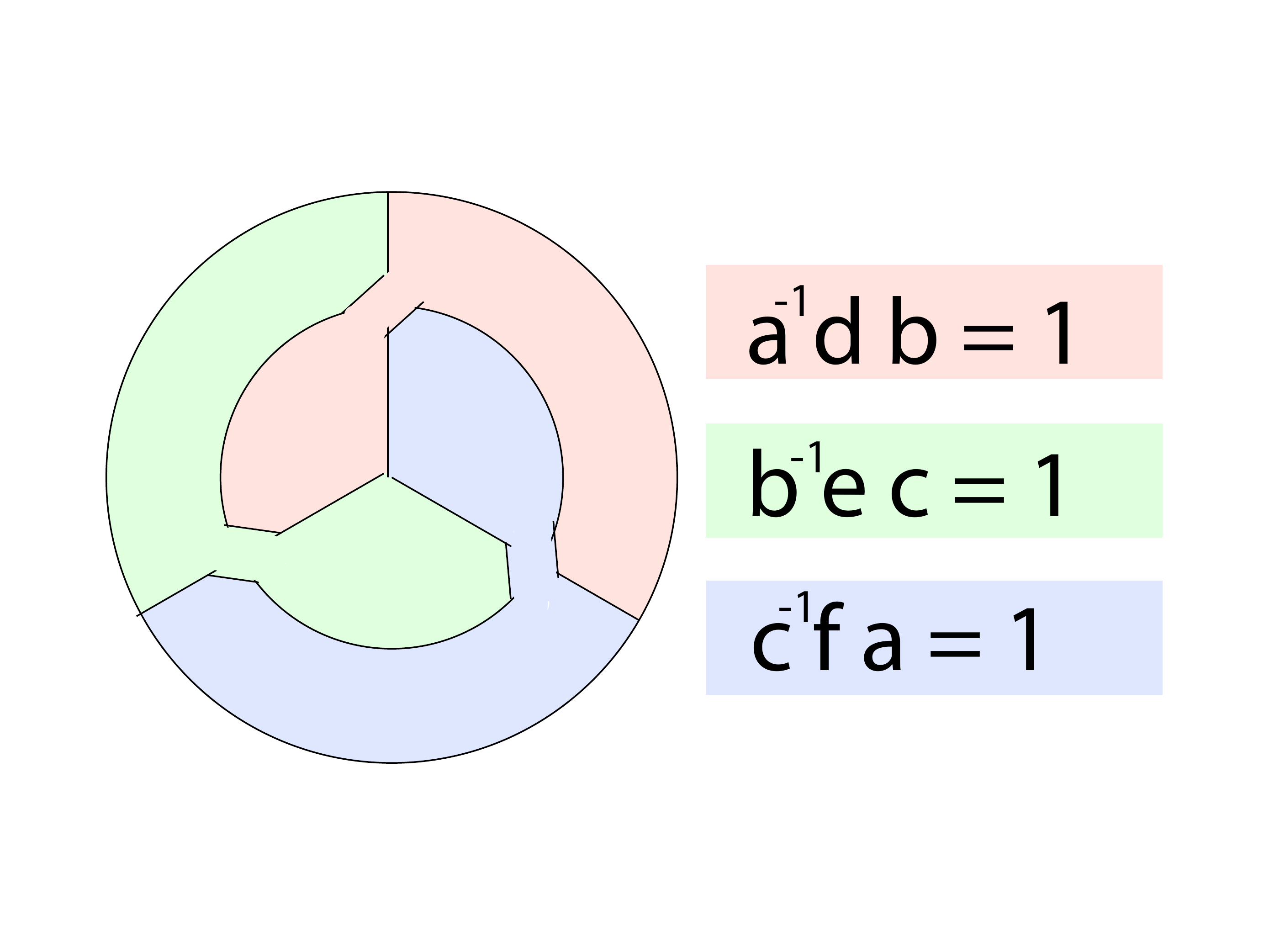}} &
\raisebox{-0.18 \columnwidth}{\includegraphics[width=0.2 \columnwidth]{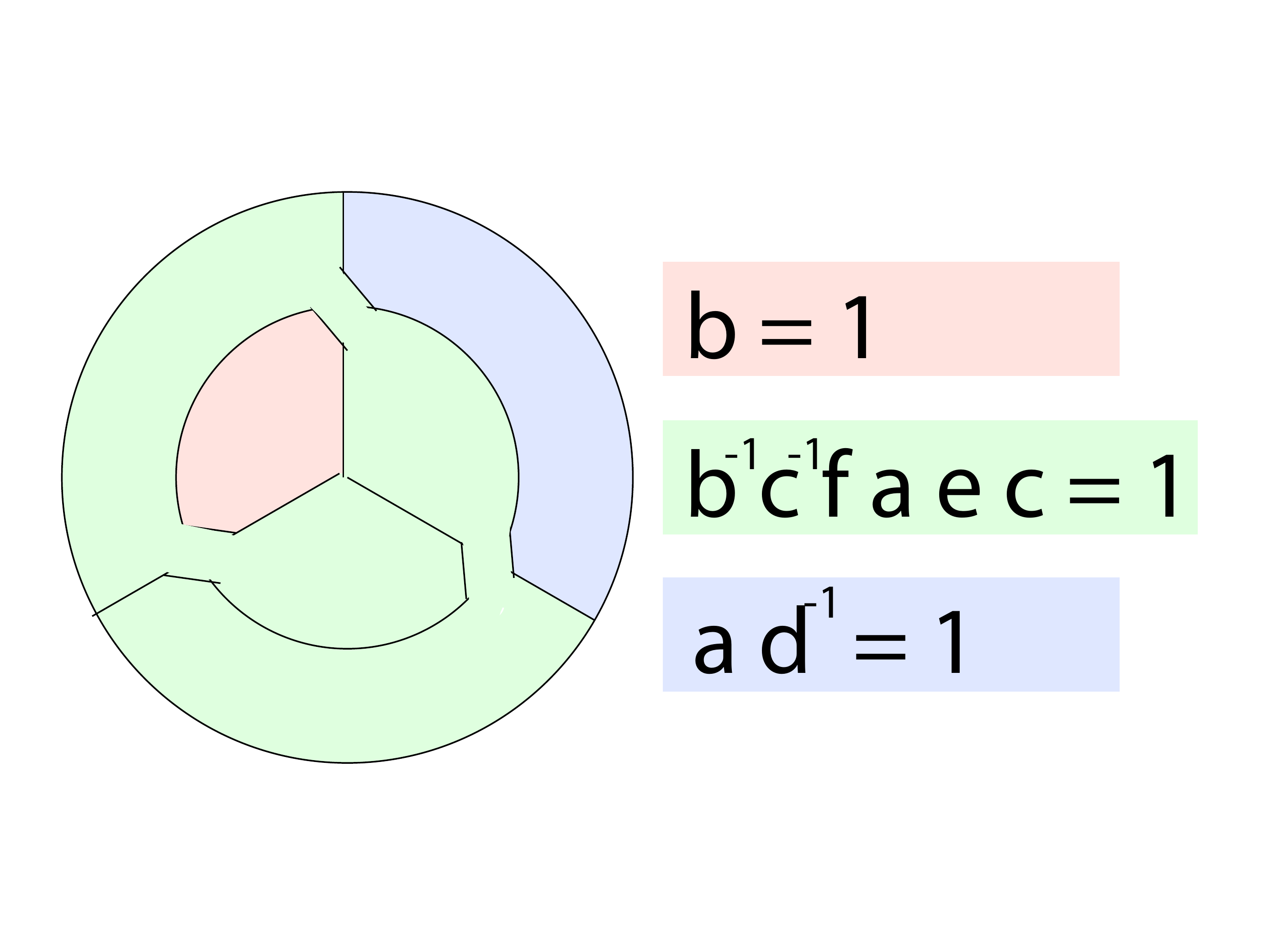}} & Step 3: Same for the lower graph. Notice that both upper and lower graphs are different in structure for the trefoil and the unknot. The trefoils' graphs have three-fold symmetry, and the facets all overlay two facets of the medial graph. The unknot's graphs have a great spread, from facets that occupy a single medial facet, to a facet that sprawls around so much as to self-abut. 
\\\hline
\parbox[t]{2 cm}{ \flushright $aeb^{-1} = 1 $
$bfc^{-1} = 1 $
$cda^{-1} = 1 $
$a^{-1}db = 1 $
$b^{-1}ec = 1 $
$c^{-1}fa = 1 $
}
&
\parbox[t]{2.4 cm}{ \flushright 
$a = 1 $\\
$b = 1 $ \\
$be^{-1} = 1 $
$ad^{-1}= 1$ 
$cdbfc^{-1}a^{-1}= 1 $
$b^{-1}c^{-1}faec = 1 $
}
&
The relations tables are then simplified. Solving for $d,e,$ and $f$ in the lower three equations for the trefoil and  substituting into the top three gives the same set of equations that the Wirtinger Presentation method yields, from which follows the group presentation for the trefoil knot, $<  b,c \, | \, bcb = cbc >$.
For the unknot, $a$ and $b$ fall out of the group presentation immediately, rapidly followed by $d$ and $e$, leaving behind the relations $cfc^{-1} = 1$ and $c^{-1}fc = 1$, which hold the same algebraic content, namely $f = 1$.  Therefore the presentation of the fundamental group in this case is $< c >$, i.e. that $\pi_1(K) \cong \mathbf{Z}$ as it should for the unknot. \\\hline
$<b,c  |  bcb = cbc>$ & 
$< c >$ & Final form \\\hline
$(3 3 3 3 3 3)^{1/6} = 3$ &
$(1 1 2 2 6 6)^{1/6} = 2.29$ &
The complexity of the presentations, defined as the \textit{geometric mean} of the relation lengths. The lower complexity for the unknot reflects the fact that, despite having the same algebraic mean in relation length, the unknot has much higher variance in length, leading to two relations that simplify immediately and cause a cascade of simplifications. 
\\
\hline
\end{tabular}

\end{table}

Finally, the mid-band set that contains the saddle surfaces that accommodate all of the knot crossings may \textit{also} be deformation retracted back onto a graph.  In this case, the full $B$ begins with sections of vertical planes above and below that match onto the ones of $C_U$ and $C_L$ respectively, but these sets of planes ``twist" into one another through the saddle-populated mid-plane of the set.  In the same way that an ``X" shape is a limiting case of hyperbolae, pulling the vertical planes back toward the mid-plane from above and below ultimately leaves a solid connection, shaped like an ``X", at each knot crossing.  $B$ is thus homotopic to the bounded medial graph of Box 1 and:

\be
\pi_1(B) \cong \mathbf{Z} \ast \cdots \ast \mathbf{Z}
\ee
where now there are twice as many copies of $\mathbf{Z}$ as in the case of the intersections, since one may view the creation of the solid connection between all four branches at a knot crossing as splitting every loop in either of the graphs associated with the intersections in two.

What are the implications of these simplifications for the inclusion-induced maps needed to form the quotient space?  Some of the sets involved in these inductions have now been shown to have trivial fundamental groups, and therefore some of the associated $i_{AB}$s will be trivial as well.  For example, for $i_{C_UB} : \pi_1(C_U \cap B) \rightarrow \pi_1(C_U) \cong 1$, all elements must be mapped to the identity, since that is all that exists in the trivial target group.  $i_{C_UB}(\omega)$ and $i_{C_LB}(\omega)$ are hence both homotopic to the null path and in the conjugacy class of the identity.

\begin{figure}[tb]
\begin{center}
\resizebox{8cm}{!}{\includegraphics{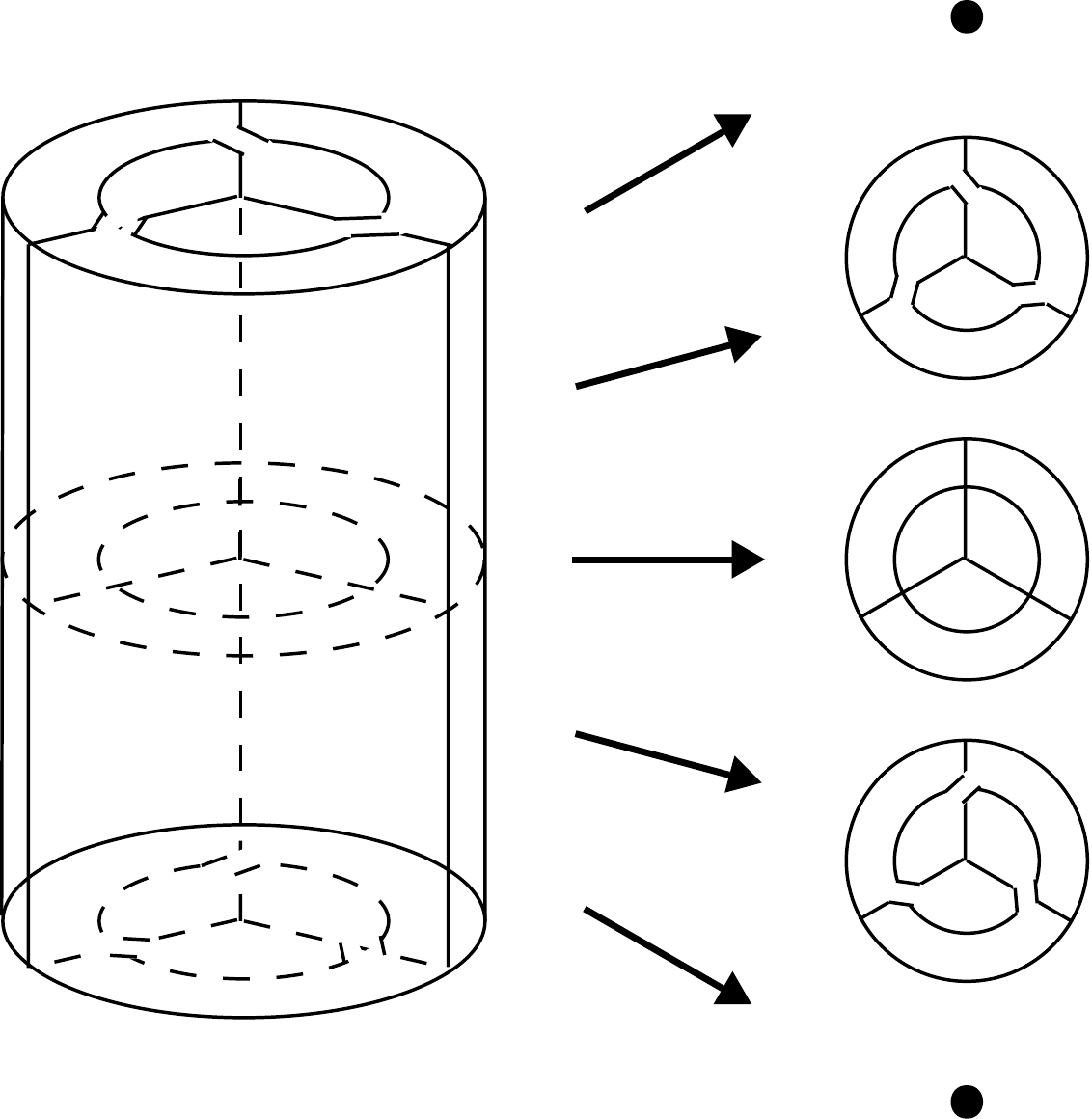}}
\end{center}
\caption{The partitioning of the template -- shown here for the trefoil knot -- into path-connected open sets amenable to application of the van Kampen Theorem.}
\label{fig:vanKampen}
\end{figure}

With these reductions in hand we may now re-examine Eqs. \ref{vK1} and \ref{vK2}:
\bea
\pi_1(K) \cong \pi_1(\textit{template}) &\cong& \pi_1(C_U \cup B) / [i_{BC_L}(\omega)] \\
\pi_1(C_U \cup B) &\cong& \pi_1(B) / [i_{BC_U}(\omega)].
\eea
or, in other words, $\pi_1(K)$ is isomorphic to the free group of $2n$ copies of $\mathbf{Z}$ modulo the generators of the free groups associated with each of the two intersection regions, where $n$ is the number of generators for each of the intersection regions' (free) fundamental groups.  Example calculations for $\pi_1(K)$ for the trefoil knot and a trefoil-like presentation of the unknot are given in Box 2.

These example calculations for the group presentation have an evident dissimilarity in the spread of the length of each relation. While the average number of symbols in all presentations remains the same (18), and therefore the average word length is constant, the unknot has much bigger variance, from which the likelihood of having relations of length 1 that simplify immediately is larger. The chances that a randomly presented group will unravel like the presentation of the unknot is related to the length of the presentations \cite{olivier}; it follows from some deep results of Gromov that this probability is governed by a logarithmically measured density of the relations in the space of all words \cite{gromov}. In our case this implies that we should measure not the algebraic mean of the word length, but the geometric mean of word lengths. Evidently, for the trefoil this number is still $3$, but for the presentation of the unknot we obtain $12^{1\over 3} = 2.289...$. This motivates us to define the \textit{presentation complexity} of a knot as the geometric mean of the number of letters in its group's presentation as defined here. Because the convexity of the logarithm implies that for a given algebraic mean the geometric mean is highest when there is a minimal variance, the presentational complexity is highest for alternating presentations.

\section*{Discussion}\label{sect:discussion}

We have shown that all the topological information carried by the knotspace, and some interesting geometry besides can be condensed into a lower dimensional object, the template.  We derived a novel, geometric presentation of the fundamental group based on the template. Conversely, the template gives us a fairly closed-form fashion in which {\em every knotspace} can be represented:  the template can be explicitly described as a set of polygons erected over the edges of the medial graph, plus saddles (in one of two choices) on the four-fold vertices, and the medial graph itself can be characterised in closed form as a bipartite planar graph with square facets and a partition of four-fold vertices. 

Owing to the geometric simplicity of the construction and the fact that the number of polygons required is linear in 
$C$, we believe that our template is the ideal tool for analyzing complicated, tangled, possibly knotted curves and links that have begun to pique interest in multiple disciplines, from polymer science to biophysics, such as linked and knotted DNA, to the study of photonic phase fields and beyond.  
Meanwhile, on the purely knot-theoretic side of the ledger, that same algorithmic simplicity and added geometric information, together with the results of several low-$C$ trials provide hope for a possibly sub-exponential-time solution to the unknotting problem. 

\begin{acknowledgements}
This work was supported in part by the National Science Foundation under grant ID~1058899, and by the Simons Foundation.
\end{acknowledgements} 

\bibliographystyle{rspublicnat}

\begin{thebibliography}{99}

\providecommand{\natexlab}[1]{#1} \expandafter\ifx\csname
urlstyle\endcsname\relax
  \providecommand{\doi}[1]{doi:\discretionary{}{}{}#1}\else
  \providecommand{\doi}{doi:\discretionary{}{}{}\begingroup
  \urlstyle{rm}\Url}\fi



\bibitem{Thurston} Thurston, W. P., \textit{The Geometry and Topology of 3-Manifolds}, Princeton Math. Dept. (1978).

\bibitem{Witten} Witten, E., Quantum field theory and the Jones polynomial. \textit{Commun. Math. Phys.} {\bf 121} 351 -- 399 (1989).

\bibitem{field theory} Faddeev, L. and Niemi, A. J.,  Stable knot-like structures in classical field theory.
\textit{Nature}  {\bf 387}.6628 58 -- 61 (1997).

\bibitem{statmech1} Frankkamenetskii, M. D., Lukashin, A. V. and Vologodskii, A. V., Statistical-mechanics and topology of polymer-chains. \textit{Nature} {\bf 258}.5534 398 -- 402 (1975).

\bibitem{Dennis} Leach, J., Dennis, M. R., Courtial, J. and Padgett, M. J., Laser beams: knotted threads of darkness. \textit{Nature} {\bf 432} 165 (2004).

\bibitem{protein} Wagner, J. R., Brunzelle, J. S., Forest, K. T. and Vierstra, R. D., A light-sensing knot revealed by the structure of the chromophore-binding domain of phytochrome. \textit{Nature} {\bf 438}.7066 325 -- 331 (2005).

\bibitem{rope length} Katritch, V.,  Bednar, J., Michoud, D., Scharein, R. G., Dubochet, J. and Stasiak, A., 
Geometry and physics of knots. \textit{Nature} {\bf 384}.6605 142 -- 145 (1996).

\bibitem{topoisomerase} Yan, J., Magnasco, M. O. and Marko, J. F., A kinetic proofreading mechanism for disentanglement of DNA by topoisomerases. \textit{Nature} {\bf 401} 932 -- 935 (1999).

\bibitem{Haken} Haken, W., Some results on surfaces in 3-manifolds, in \textit{Studies in Modern Topology.} P. J. Hilton, ed, Math. Assoc. Amer., Prentice-Hall, Englewood Cliffs, NJ, 39 -- 98 (1968).  

\bibitem{classification1} Hemion, G., \textit{The classification of knots and 3-dimensional spaces}. Oxford Univ. Press, Oxford (1992). 

\bibitem{classification2} Hoste, J., The enumeration and classification of knots and
links, in \textit{Handbook of Knot Theory}. W. W. Menasco and M. B. Thistlethwaite, eds, Elsevier B. V., Amsterdam (2005). 

\bibitem{classification3} Dowker, C. H. and Thistlethwaite, M. B., On the classification of knots. \textit{C. R. Math. Rep. Acad. Sci. Canada} {\bf 4} no. 2, 129 -- 131 (1982).

\bibitem{tabulation} Rankin, S., Flint, O. and Schermann, J., Enumerating the prime alternating knot. 
\textit{J. Knot Theory Ramificat.} {\bf 13}(1), part I: 57 -- 100, part II 101 -- 150 (2004).

\bibitem{spines} Matveev, S. V., \textit{Algorithmic Topology and Classification of 3-Manifolds.} Springer-Verlag, Berlin (2007).

\bibitem{braids} Birman, J. S. and Hirsch, M. D., A new algorithm for recognizing the unknot. \textit{ Geometry and Topology} {\bf 2}(9) 175 -- 220 (1998).

\bibitem{stateart2} Liu, Z. and Chan H. S., Efficient chain moves for Monte Carlo simulations of a wormlike DNA model: Excluded volume, supercoils, site juxtapositions, knots, and comparisons with random-flight and lattice models. \textit{J. Chem. Phys.} {\bf 128} 145104 (2008).


\bibitem{Wirtinger} Wirtinger, W., \"{U}ber die Verzweigungen bei Funktionen von zwei Verl\"{a}nderlichen, \textit{Jahresbericht d. Deutschen Mathematiker Vereinigung.} {\bf 14} 517 (1905).

\bibitem{Hatcher} Hatcher, A., \textit{Algebraic Topology} Cambridge Univ. Press, Cambridge (2002).

\bibitem{Floer1} Ozsvath, P. and Szabo, Z., Holomorphic disks and knot invariants. \textit{Adv. Math.}\textbf{186}(1) 58 (2004).

\bibitem{Floer2} Rasmussen, J., Floer homology and knot complements. arXiv:math/0306378 [math.GT] (2003).

\bibitem{Kirkwood} Kirkwood, J. G., Critique of the free volume theory of the liquid state. \textit{J. Chem. Phys.} \textbf{18} 380 (1950).

\bibitem{stat} Salsburg, Z. W. and Wood, W. W., Equation of state of classical hard spheres at high density. \textit{J. Chem. Phys.} \textbf{37} 798 (1962).

\bibitem{Voronoi} Okabe, A., Boots, B., Sugihara, K. and Chiu, S. N., \textit{Spatial Tesselations -- Concepts and Applications of Voronoi Diagrams}, 2nd Ed. John Wiley \& Sons, Chichester (2000).

\bibitem{medial} Blum, H., A transformation for extracting new descriptors of shape, in \textit{Models for the perception of speech and visual form.} M.I.T. Press, Cambridge, MA, 362 (1967).

\bibitem{pi1prime} Simon J., How many knots have the same group? \textit{Proc. Amer. Math. Soc.} {\bf 80}(1) 162 (1980).

\bibitem{olivier} Ollivier, Y., \textit{A January 2005 invitation to random groups}, Sociedade Brasileira de Matematica (2005).

\bibitem{gromov} Gromov, M.,  Asymptotic invariants of infinite groups, in \textit{Geometric group theory}, G. Niblo and M. Roller, eds, Cambridge University Press, Cambridge (1993).

\end{thebibliography}

\end{document}